\documentclass[ejs]{imsart}

\usepackage{amsfonts}
\usepackage{amsthm,amsmath}
\usepackage{bbm, mathtools}
\usepackage{enumerate}
\usepackage{caption}
\usepackage{here}
\usepackage{color}

\doi{10.1214/154957804100000000}
\pubyear{0000}
\volume{0}
\firstpage{1}
\lastpage{1}
\startlocaldefs

\def\R{{\mathbb{R}}}

\def\N{{\mathbb{N}}}

\def\P{{\mathbb{P}}}

\newcommand{\red}[1]{\color{red}   #1 \color{black} }

\DeclarePairedDelimiter\norm{\lvert\lvert}{\rvert\rvert}
\DeclareMathOperator*{\argmin}{arg\,min}

\newtheorem{Lemma}{Lemma}[section] %numbering by section
\newtheorem{Theorem}{Theorem}[section] %numbering by section
\newtheorem{Assumption}{Assumption}[section] %numbering by section
\newcommand{\proofend}{$\hfill\Box{~}$}
\newenvironment{Proof}{\noindent {\em{\bf Proof.\ }}}{\proofend\\}

\endlocaldefs

\begin{document}

\begin{frontmatter}

% "Title of the Paper"
\title{Robust censored regression with $\boldsymbol{\ell}_1$-norm regularization}
\runtitle{Robust censored regression}

% indicate corresponding author with \corref{}
% \author{\fnms{John} \snm{Smith}\thanksref{t1}\corref{}\ead[label=e1]{smith@foo.com}\ead[label=e2,url]{www.foo.com}}
% \thankstext{t1}{Thanks to somebody} 
% \address{line 1\\ line 2\\ \printead{e1}\\ \printead{e2}}

\author{\fnms{Jad} \snm{Beyhum}\ead[label=e1]{jad.beyhum@kuleuven.be}}
\and
\author{\fnms{Ingrid} \snm{Van Keilegom}\ead[label=e2]{ingrid.vankeilegom@kuleuven.be}}
\address{Research Centre for Operations Research and Statistics, KU Leuven,
Naamsestraat 69, B-3000 Leuven, Belgium\\ 
\printead{e1,e2}}
%\and
%\author{\fnms{???} \snm{???}\ead[label=e2]{ingrid.vankeilegom@kuleuven.be}}
%\address{\printead{e2}}

\runauthor{Beyhum \& Van Keilegom}

\begin{abstract}
This paper considers inference in a linear regression model with random right censoring and outliers. The number of outliers can grow with the sample size while their proportion goes to zero. The model is semiparametric and we make only very mild assumptions on the distribution of the error term, contrary to most other existing approaches in the literature. We propose to penalize the  estimator proposed by Stute for censored linear regression by the $\ell_1$-norm. We derive rates of convergence and establish asymptotic normality of the estimator of the regression coefficients. Our estimator has the same asymptotic variance as Stute's estimator in the censored linear model without outliers. Hence, there is no loss of efficiency as a result of robustness. Tests and confidence sets can therefore rely on the theory developed by Stute. The outlined procedure is also computationally advantageous, since it amounts to solving  a convex optimization program. We also propose a second estimator which uses the proposed penalized Stute estimator as a first step to detect outliers. It has similar theoretical properties but better performance in finite samples as assessed by simulations.
\end{abstract}



% history:
% \received{\smonth{1} \syear{0000}}

%\tableofcontents

\end{frontmatter}

\section{Introduction}

The present paper considers a linear model with outliers and censoring. We have at hand a dataset of $n$ independent realizations of an outcome random variable $Y_i$ and a random vector of covariates $X_i$ with support in $\mathbb{R}^{p}$, where $p$ is fixed. The variable $Y_i$ is the minimum between a random variable $T_i$ and a censoring time $C_i$, and $T_i$ and $X_i$ are related through the following model:
\begin{equation}\label{model}\begin{array}{cc}T_i=X_i^{\top}\beta+\alpha_i +\xi_i, & \quad i=1,\dots,n, \end{array}\end{equation}
where $\beta\in \mathbb{R}^p$, while the error term $\xi_i$ and $\alpha_i$ are real-valued random variables. The observation $i$ is called an outlier if $\alpha_i\ne 0$. Let the average proportion of outliers $\mathbb{P}(\alpha_i\ne 0)$ be denoted by $\pi_\alpha$. We assume that $\mathbb{E}[X_iX_i^\top]$ is positive definite and $\mathbb{E}[X_i\xi_i ]=0$. Under these conditions, $\beta =E[X_iX_i^\top]^{-1}E[X_i(T_i-\alpha_i)]$ is the coefficient of the best linear predictor of $T_i-\alpha_i$ by $X_i$. The goal is to estimate $\beta$ based on a sample $\{Y_i,X_i,\delta_i\}_i$, where  $\delta_i=\mathbbm{1}_{T_i\le C_i}$ is the censoring indicator. The vector $\alpha=(\alpha_1,\dots,\alpha_n)$ acts as a nuisance parameter that may represent measurement errors.

When the variable $T_i$ is the logarithm of a duration, model (\ref{model}) is the so-called accelerated failure time model (henceforth AFT). Applying a standard estimator for this type of models (see e.g. \cite{buckley1979linear}, \cite{koul1981regression}, \cite{stute1993consistent}, \cite{tsiatis1990estimating},\cite{akritas1996}, \cite{vkakritas2000}) would yield (under some additional assumptions) a consistent estimate of the best linear predictor of $T_i$ on $X_i$. In practice, however, because of outliers, the estimated coefficients may yield poor predictions for most observations. Robust methods seek to estimate features of the majority of the data rather than average properties of the full dataset. 

This paper studies a lasso-type estimator of the parameter vector $\beta$ that is robust to outliers. Our proposal penalizes the criterion of the inverse probability weighting estimator defined in \cite{stute1993consistent} by the $\ell_1$-norm. In an asymptotic regime where the average proportion of outliers $\pi_\alpha$ goes to $0$ with $n$ while the mean number of outliers $n\pi_\alpha$ is allowed to diverge to $\infty$, we derive rates of convergence of the estimator.  If $\pi_\alpha$ goes to $0$ quickly enough the suggested estimator is asymptotically normal with the same asymptotic variance as Stute's estimator in the model without outliers. Hence, no price is paid for robustness. We also develop a two-step procedure, which uses the $\ell_1$-norm penalized estimator to detect outliers, and then run Stute's  estimator on the sample that is cleaned of outliers. We present a simple computational algorithm for our lasso-type estimator. Our methods exhibit good finite sample properties in simulations. The proposed approaches are simple, and more importantly, do not rely on a parametric assumption on the error term $\xi_i$.\\

\noindent\textbf{Related literature.} Although many papers have studied robust estimation in the Cox model (see e.g. \cite{sasieni1993maximum, sasieni1993some, bednarski1999adaptive, bednarski2003robustness, bednarski1993robust, lin1989robust, farcomeni2011robust}), far fewer papers have considered the AFT model with outliers. In \cite{locatelli2011robust} a three-step procedure is proposed that has a high breakdown point and is efficient when there are no outliers. A preliminary high breakdown point $S$-estimator is used to detect and remove outliers from the dataset. In contrast, \cite{sinha2019robust} studies an approach that bounds the influence function in both the outcome and the covariates. It is robust in both the $x$ and the $y$ dimensions. Both methods require a parametric assumption on the distribution of the error term $\xi_i$. Instead, our model is semiparametric. Finally, in a semiparametric model, \cite{heller2007smoothed} proposes an approach based on smooth weighted rank estimation. Unlike in our paper, the error term is assumed to be independent of the regressors. The influence function of the estimator is bounded. However, the estimator is less efficient than non-robust rank-based estimation techniques, that is, unlike with our $\ell_1$-regularization approach,  some price is paid for robustness.

Another related field is that of robust estimation with lasso-type estimators. Many papers have studied a $\ell_1$-penalized least squares estimator for (uncensored) contaminated linear regression (see e.g. \cite{gannaz2007robust, she2011outlier, lee2012regularization,lambert2011robust,dalalyan2012socp,gao2016penalized, collier2017rate}). In this setup, \cite{beyhum2020inference} develops estimation results similar to those of the present paper but does not allow for censoring and does not study the two-step procedure.\\

\noindent\textbf{Outline.} In Section \ref{sec.est}, we present the estimators, the assumptions and the convergence results.  An algorithm to compute the estimators and the results of a simulation study are discussed in Section \ref{sec.sim}. Finally, all technical details and proofs are given in the Appendix. \\

\noindent\textbf{Notation.} We use the following notation. For a matrix $M$, $M^{\top}$ is its transpose, $\norm{M}_2$,  $\norm{M}_1$ and $\norm{M}_{\infty}$ are the $\ell_2$-norm, $\ell_1$-norm and the sup-norm of the vectorization of $M$, respectively, $\norm{M}_{\text{op}}$ is the operator norm of $M$ and $\norm{M}_0$ is the number of non-zero coefficients in $M$, that is its $\ell_0$-norm. Moreover, $M_{k\cdot}$ denotes its $k^{th}$ column. If there is a second matrix $A$ of the same size as $M$, then $\left<M,A\right>$ is the Frobenius scalar product of $M$ and $A$. For a set $S\subset \{1,\dots,p\}$ and a vector $v\in\R^p$, let $|S|$ denote the number of elements in $S$ and $v_{S}\in \R^p$ the vector such that $(v_S)_i=v_i$ for all $i\in S$ and $(v_S)_i=0$ otherwise. We also introduce $v_{S^c}=v-v_S$.\\

\section{Estimation results}\label{sec.est}
\subsection{Probabilistic framework} 
We consider a sequence of data generating processes (henceforth, DGPs) depending on the sample size $n$. The joint distribution of $(X_i,\xi_i)$ does not depend on $n$, but the distribution of $\alpha_i$, and therefore also of $Y_i$, do depend on $n$.
%is fixed with $n$. 
We study a regime in which  $n$ goes to $\infty$ and the contamination level $\pi_{\alpha}$ goes to $0$  with $n$. This implies that asymptotically only a negligible proportion of observations is not generated by the linear model $T_i=X_i^\top \beta+\xi_i$. For $i=1,\dots n$, let $\tilde T_i=T_i-\alpha_i$, $\tilde Y_i=\min(\tilde T_i,C_i)$ and $\tilde \delta_i =1_{\{\tilde T_i\le C_i\}}$. This is the sample that would be observed if there were no outliers. Since $\pi_\alpha\to 0$, the proportion of observations from this conceptual sample which are part of the observed sample goes to $1$.

\subsection{Main estimator}
\label{sec.estim1}
Note that in this asymptotic setting, even if there were no censoring, the ordinary least squares (OLS) estimator would not be consistent. Indeed, if we assume that $X_i$ is a constant equal to $1$ in model \eqref{model}, then, conditional on the $\alpha_i$'s, the average value of the OLS estimator of the regression of $T_i$ on the constant would be  $\beta+(\sum_{i=1}^n \alpha_i)/n$, which would diverge if the nonzero coefficients of $\alpha$ are too large. Censoring poses additional challenges because $T_i$ is not observed for some $i$. It is well known that $\ell_1$-norm penalization allows to build robust estimators (\cite{she2011outlier,beyhum2020inference}). In \cite{stute1993consistent}, Stute proposes a simple estimator for censored regression. In this paper, we propose to combine these two ideas: our estimator penalizes the criterion of the Stute estimator by the weighted sum of the absolute values of the coefficients of $\alpha$.

Let $Y_{(1)},\dots, Y_{(n)}$ denote the order statistics of $\{Y_i\}_{i=1}^n$. For a vector $R\in\R^n$, $R_{(i)}$ is the coordinate of $R$ associated with $Y_{(i)}$. For the sake of simplicity, we assume that the distribution of $Y_i$ is continuous. Let us introduce the Kaplan-Meier weights:
\begin{equation} \label{kmweights} w_{(1)}= \frac{\delta_{(1)}}{n},\ w_{(i)}= \frac{\delta_{(i)}}{n-i+1}\prod_{j=1}^{i-1}\left( \frac{n-j}{n-j+1}\right)^{\delta_{(j)}},\ i=1,\dots,n,
\end{equation}
and the estimator:
\begin{equation} \label{estimator}
(\widehat{\beta},\widehat{\alpha})\in\argmin\limits_{b\in\R^p,\ a\in\R^n}\sum_{i=1}^n w_{(i)} (Y_{(i)} -X_{(i)}^\top b- a_{(i)})^2+ \lambda \sum_{i=1}^n\sqrt{w}_{(i)}|a_{(i)}|,
\end{equation}
where $\lambda>0$ is a penalty level.
Remark that the penalty that we use weights the coefficients of the vector $a$ by the square-root of the Kaplan-Meier weights, this is because the latter determine the influence of the entries of $a$ on $\sum_{i=1}^n w_{(i)} (Y_{(i)} -X_{(i)}^\top b- a_{(i)})^2$. In particular, when $w_{(i)}=0$, then $a_{(i)}$ should not be penalized. Let us introduce further notations: $X^w=(\sqrt{w}_{(1)}X_{(1)},\dots, \sqrt{w}_{(n)}X_{(n)})^\top$, $w=(w_{(1)},\dots, w_{(n)})^\top$ and for $a\in\R^n$, define $a^w=(\sqrt{w}_{(1)}a_{(1)},\dots, \sqrt{w}_{(n)}a_{(n)})^\top$. The estimator \eqref{estimator} can be rewritten as
\begin{equation} \label{estimator2}
(\widehat{\beta},\widehat{\alpha})\in\argmin\limits_{b\in\R^p,\ a\in\R^n}\left|\left|Y^w-X^wb-a^w \right|\right|_2^2+ \lambda \norm{a^w}_1.
\end{equation}
For $a\in\R^n$, let $\widehat{\beta}(a)$ be the weighted least squares estimator 
\begin{equation} \label{estimatorprofile}
\widehat{\beta}(a)\in\argmin\limits_{b\in\R^p} \left|\left|Y^w-X^wb-a^w \right|\right|_2^2.
\end{equation}
An important remark is that
\begin{equation}\label{profile}\left|\left|Y^w-X^w\widehat{\beta}-\widehat{\alpha}^w \right|\right|_2^2+ \lambda \norm{\widehat{\alpha}^w}_1\le \left|\left|Y^w-X^wb-\widehat{\alpha}^w \right|\right|_2^2+ \lambda \norm{\widehat{\alpha}^w}_{1},\end{equation}
for any $b\in \mathbb{R}^{p}$ and, therefore, $\widehat{\beta}=\widehat{\beta}(\widehat{\alpha})$ when there is a unique solution to the minimization program \eqref{estimatorprofile}. Hence, when $(X^w)^\top X^w$ is positive definite, we have
\begin{equation}\label{wOLS}\widehat{\beta}=\left((X^w)^\top X^w\right)^{-1}(X^w)^\top (Y^w-\widehat{\alpha}^w).\end{equation}
Then, let $P_{X^w}=X^w\left((X^w)^\top X^w\right)^{-1}(X^w)^\top$ be the orthogonal projector on the columns of $X^w$ and let $M_{X^w}=I_n-P_{X^w}$. For all $a\in \mathbb{R}^n$ and $b \in \mathbb{R}^{p}$, we have 
$$\left|\left| M_{X^w}(Y^w-a^w)\right|\right|_2^2+\lambda\norm{a^w}_1\le\left|\left| Y^w -X^wb -a^w\right|\right|_2^2+\lambda\norm{a^w}_1.$$
Therefore, since $\left|\left| Y^w -X^wb -a^w\right|\right|_2^2=\left|\left| M_{X^w}(Y^w-a^w)\right|\right|_2^2$ if $X^wb=P_{X^w}(Y^w-a^w)$, it holds that
\begin{equation} \label{progalpha} \widehat{\alpha}\in \argmin_{a\in \mathbb{R}^n}\left|\left| M_{X^w}(Y^w-a^w)\right|\right|_2^2+\lambda\norm{a^w}_1.\end{equation}
The expressions (\ref{wOLS}) and (\ref{progalpha}) are useful in the proofs of our theoretical results.  However, for the computation of the estimator we will use an algorithm that is outlined in Section \ref{sec.algo}.

\subsection{Assumptions}
In this section, we state the different assumptions that we need to prove the asymptotic normality of our estimator $\widehat\beta$.  The first set of assumptions is standard in linear models and would ensure asymptotic normality of the Stute estimator of the regression of $\tilde T_i$ on $X_i$.

\begin{Assumption}\label{DGP} The following holds:
\begin{enumerate}
 \item[(i)]$\{C_i,X_i,\xi_i\}_i$ are independent random variables and for fixed $n$ they have the same distribution;
 \item[(ii)]$ \mathbb{E}[X_i\xi_i]=0$;
 \item[(iii)]$\Sigma_X=\mathbb{E}[X_iX_i^\top]$ exists and is positive definite;
\item[(iv)] for all $t>0$, $\norm{\xi}_\infty=o_P(n^t)$.
 \end{enumerate}
\end{Assumption}
Assumption \ref{DGP}(iv) is a mild condition on the tails of the distribution of the error term. Many usual distributions of the accelerated failure time model (Gaussian, Laplace, logistic, Weibull and gamma) satisfy it (this can be shown using tail bounds).

The next assumption concerns the choice of the tuning parameter.
\begin{Assumption}\label{effective}
Let $\lambda =n^{\lambda_0-(\widehat{\pi}_{uc}/2)}$, where $\lambda_0>0$ and $\widehat{\pi}_{uc} =\sum_{i=1}^n \delta_i/n$.
\end{Assumption}
The quantity $\widehat{\pi}_{uc}$ is an estimator of the probability that an observation is uncensored, denoted $\pi_{uc}=\P(\tilde \delta_i=1)$. Note that $\lambda_0$ can be chosen arbitrarily small. In the simulations, we will set $\lambda_0=10^{-4}$.

Next, we have an assumption on the censoring mechanism. For a random variable $R$, let $\tau_R$ denote the upper bound of the support of its distribution.
\begin{Assumption}\label{censoring} The following holds:
\begin{enumerate}
 \item[(i)] The censoring variable $C_i$ is independent of $(T_i,\alpha_i,X_i)$;
\item[(ii)] $\tau_{C_1}=\infty$ or $\tau_{\tilde T_1}< \tau_{C_1}$.
 \end{enumerate}
\end{Assumption}
The first condition is the usual independent censoring assumption from the survival analysis literature, while the second condition is a sufficient follow-up assumption. They can be relaxed, for instance  \cite{stute1993consistent} and \cite{stute1996distributional} make weaker (but less understandable) assumptions. 

Moreover, we assume that some rate conditions from \cite{stute1996distributional}  hold. These conditions are discussed in \cite{stute1995central} and \cite{stute1996distributional} and correspond to some of the requirements for asymptotic normality of Stute's estimator in the absence of outliers.. We state them in Assumption \ref{stute_rate} in Appendix \ref{app.cond}.

Finally, we also make the following additional rate conditions on $\norm{X}_\infty, \lambda, \pi_\alpha$ and $\pi_{uc}$:
\begin{enumerate}[(A)]
 \item\label{Ei} $n^{(1-\pi_{uc}+3\lambda_0)/2} \norm{X}_\infty \sqrt{\pi_\alpha} =o_P(1)$;
 \item\label{Eii} $n^{(3/2)-\pi_{uc}+3\lambda_0} \norm{X}_\infty \pi_\alpha=o_P(1)$
 \end{enumerate}
Condition (\ref{Ei}) is used to derive the rate of convergence of the estimator. We leverage (\ref{Eii}) to show asymptotic normality of the estimator. Remark that these conditions imply that $\pi_\alpha$ tends to zero, and that these rate conditions become more stringent as the censoring rate increases ($\pi_{uc}$ decreases). Therefore, the more observations are censored, the fewer outliers are allowed. Because of this, we may expect  the performance of our estimators to decrease when there is more censoring, which is what we observe in the simulations. Note that this would happen even if there were no outliers. Let us illustrate the strength of these conditions through two examples:\\

\noindent\textbf{Example 1.} Let us assume that $X_i$ are bounded (so that $\norm{X}_{\infty}=O_P(1)$, where $X=(X_1,\dots, X_n)^\top$) and $\pi_{uc}=7/8$. Then, condition (\ref{Ei}) holds when $\pi_\alpha = o(n^{-(1/8)-3\lambda_0}).$ Condition
(\ref{Eii}) is satisfied if
$\pi_\alpha = o(n^{-(5/8)-3\lambda_0})$
which allows the average number of outliers $n\pi_{\alpha}$ to go to infinity. \\

\noindent\textbf{Example 2.} Let us assume that $X_i$ are sub-Gaussian (so that $\norm{X}_{\infty} \linebreak =O_P(\sqrt{\log(n)})$) and $\pi_{uc}=3/4$. Condition (\ref{Ei}) is satisfied if $\pi_\alpha = o(n^{-(1/4)-3\lambda_0}/ \linebreak \sqrt{\log(n)}).$ Condition (\ref{Eii}) holds when
$\pi_\alpha = o(n^{-(3/4)-3\lambda_0}/\sqrt{\log(n)})$
which also allows the average number of outliers $n\pi_{\alpha}$ to go to infinity.

\subsection{Convergence results}

The following result characterizes the rate of convergence of the estimator.
 
\begin{Theorem} \label{consth} Let Assumptions \ref{DGP}, \ref{effective}, \ref{censoring}, \ref{stute_rate}  and the rate condition (\ref{Ei}) hold. Then, we have: 
$$\norm{\widehat{\beta}-\beta}_2=O_P\left(\max\left(\frac{1}{\sqrt{n}}, n^{1-\pi_{uc}+3\lambda_0}\norm{X}_\infty\pi_\alpha\right)\right).$$
\end{Theorem}

Under Example 1, the rate becomes $\max(n^{-\frac12}, n^{(1/8)+3\lambda_0}\pi_\alpha)$. The next theorem states that the estimator is asymptotically normal.

\begin{Theorem} \label{an} Let Assumptions \ref{DGP}, \ref{effective}, \ref{censoring}, \ref{stute_rate}  and rate conditions (\ref{Ei}) and (\ref{Eii}) hold. Then, we have: 
$$\sqrt{n} (\widehat{\beta}-\beta )\xrightarrow{d}  \mathcal{N}(0,   \Sigma_X^{-1}\Sigma\Sigma_X^{-1}),$$ where $\Sigma$ is defined in Appendix \ref{app.as}.

\end{Theorem}
Note that the asymptotic variance of our estimator is the same as that of \cite{stute1996distributional} for the regression of $T_i-\alpha_i$ (censored by $C_i$) on $X_i$. Hence, there is no loss in efficiency, while our estimator remains asymptotically normal in the presence of outliers. This theorem allows to build confidence intervals and tests on $\beta$.
These confidence intervals are obtained under an asymptotic regime with triangular array data where the number of outliers is allowed to go to infinity while their proportion goes to $0$. A 95\% confidence interval I on a functional $L(\beta)$ of $\beta$ built with Theorem \ref{an} should therefore be interpreted as follows: if the number of outliers in our data is low enough and the sample size is large enough,
then there is a probability of approximatively 0.95 that $L(\beta)$ belongs to I.

\subsection{Second step estimator}

Let us now consider the following second step estimator. Pick a threshold level $\tau_0\ge 0$ and define $$\widehat{J}(\tau_0)=\{i\in\{1,\dots,n\}:\ |\widehat{\alpha}_{(i)}^w|> \tau_0\},$$ the set of outliers detected by $\widehat{\alpha}$ at level $\tau_0$. The role of the threshold is to allow for small mistakes of the first step estimator.  In practice it should be chosen so that very small coefficients of $\widehat{\alpha}^w$ are not considered as outliers.
The second-step estimator corresponds to the weighted by $w$ regression of $Y$ on $X$  only on the observations that do not belong to $\widehat{J}(\tau_0)$, that is 
$$\widetilde{\beta}\in \argmin\limits_{b\in\R^p}\sum_{i\notin \widehat{J}(\tau_0)}^n w_{(i)} (Y_{(i)} -X_{(i)}^\top b)^2.$$
We can define an estimator $\widetilde{\alpha}$ such that $$(\widetilde{\beta},\widetilde{\alpha})\in \argmin\limits_{ b\in\R^p,a\in\R^n}\norm{Y^w-X^wb - a_{\widehat{J}(\tau_0)}^w}_2^2, $$
and we have $\widetilde{\alpha}^w_{\widehat{J}(\tau_0)}=\big(Y^w-X^w\widetilde{\beta}\big)_{\widehat{J}(\tau_0)}$. Remark also that by arguments similar to that of the end of Section \ref{sec.estim1}, it holds that
\begin{equation}\label{equiv_formulation} \widetilde{\beta} \in \argmin\limits_{b\in\R^p}\norm{Y^w-X^wb -\widetilde{\alpha}_{\widehat{J}(\tau_0)}^w}_2^2;\ \widetilde{\alpha}\in \argmin_{a\in \mathbb{R}^n}\left|\left| M_{X^w}(Y^w- a_{\widehat{J}(\tau_0)}^w)\right|\right|_2^2.\end{equation}
The following theorem states that the one-step and two-step estimators are asymptotically equivalent.

\begin{Theorem} \label{an2} Let Assumptions \ref{DGP}, \ref{effective}, \ref{censoring}, \ref{stute_rate} and condition (\ref{Ei}) hold. Then, we have:
$$\norm{\widetilde{\beta}-\widehat{\beta}}_2=O_P\left(n^{1-\pi_{uc}+3\lambda_0}\norm{X}_\infty\pi_\alpha\right),$$ which implies that $\widetilde{\beta}$ satisfies the rate of convergence given in Theorem \ref{consth}. Moreover, if condition (\ref{Eii}) also holds, we have: 
$$\sqrt{n} (\widetilde{\beta}-\widehat{\beta})=o_P\left(1\right),$$
which yields that $\widetilde{\beta}$ follows Theorem \ref{an} (with $\widehat{\alpha}$ replaced by $\widetilde{\alpha}_{\widehat{J}(\tau_0)}$).
\end{Theorem}

As we will see the main advantage of the two-step estimator is that it works better in simulations. This is because $\widetilde{\alpha}$ does not suffer from the shrinkage bias imposed on $\widehat{\alpha}$ by the $\ell_1$ penalty.

\section{Computation and simulations}\label{sec.sim}
\subsection{Iterative algorithm}\label{sec.algo}

We propose to use an iterative algorithm over $b,a$ to compute our estimator \eqref{estimator2}. Let us start from $\left(b^{(0)},a^{(0)}\right)$ and compute the following sequence for $t\in\mathbb{N}^*$ until convergence:
\begin{enumerate}
\item $b^{(t+1)}\in \argmin_{b\in \mathbb{R}^{p}} \norm{ Y^w-X^wb-(a^{(t)})^w}_2^2;$
\item $a^{(t+1)}\in \argmin_{a\in \mathbb{R}^n}  \norm{ Y^w-X^wb^{(t+1)}-a^w}_2^2+\lambda\norm{a^w}_1.$
\end{enumerate}
The first step is the computation of the OLS estimator of the regression of $Y^w-(a^{(t)})^w$ on $X^w$. The following lemma is a direct consequence of Section 4.2.2 in \cite{giraud2014introduction} and shows how to compute the result of step 2.

\begin{Lemma}\label{alphaiter}
For $i=1,\dots,n$, if $\big|Y_{(i)}^w-X^w_{(i)}b^{(t+1)} \big|\le \lambda/2$ then $(a^{(t+1)})^w_{(i)}=0$. Otherwise, we have $(a^{(t+1)})^w_{(i)}=Y_{(i)}^w-X^w_{(i)}b^{(t+1)} -\text{sign}\big(Y_{(i)}^w-X^w_{(i)}b^{(t+1)}\big)(\lambda/2)$.
\end{Lemma}

\subsection{Simulations}
We consider the following DGP. There are two variables, a constant $X_{i1}=1$ and $X_{i2}$ which is uniformly distributed on the interval $[0,1]$. The error term $\xi_i$ follows a $\mathcal{N}(0,1)$ distribution. The variable $\alpha_i$ is equal to $0$ when $X_{i2}<1-(5.10^{-3})$ and $-20$ otherwise. This implies that $\pi_\alpha=5.10^{-3}$.  Since we set the sample size $n$ to $1,000$, there are on average $5$ outliers in each sample. The parameter is $\beta=(1,1)^\top$ so that $Y_i=X_{i1}+X_{i2}+\alpha_i+\xi_i$. The censoring duration $C_i$ follows a $\mathcal{N}(\mu,1)$ distribution and we vary $\mu$ in order to change the probability that an observation is not censored, namely $\pi_{uc}$. 

\begin{figure}[H]
\begin{minipage}{0.46\textwidth}
  \centering
  \includegraphics[width=65mm]{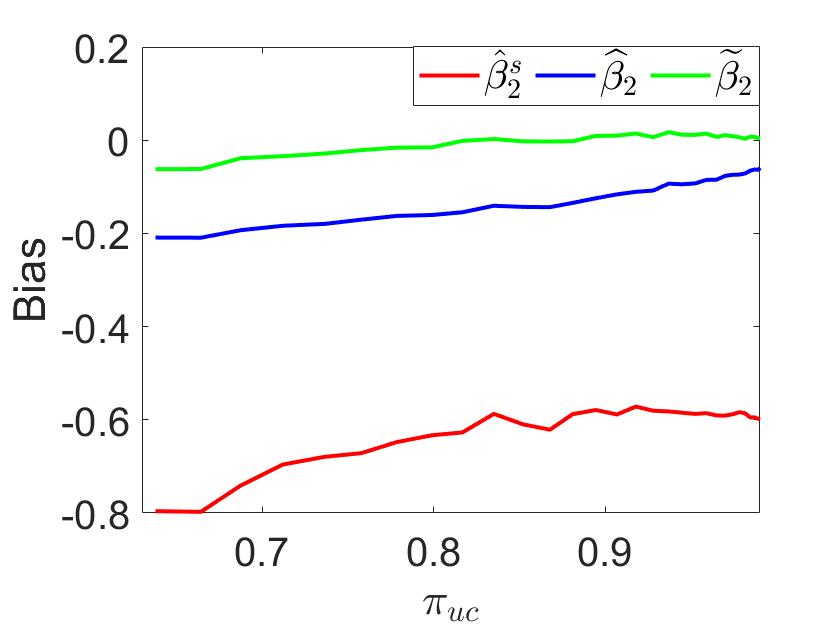}
  \caption{Bias of the estimators.}
   \label{fig:bias} 
    \end{minipage}
    \quad\quad
    \begin{minipage}{0.46\textwidth}
        \centering
   \includegraphics[width=65mm]{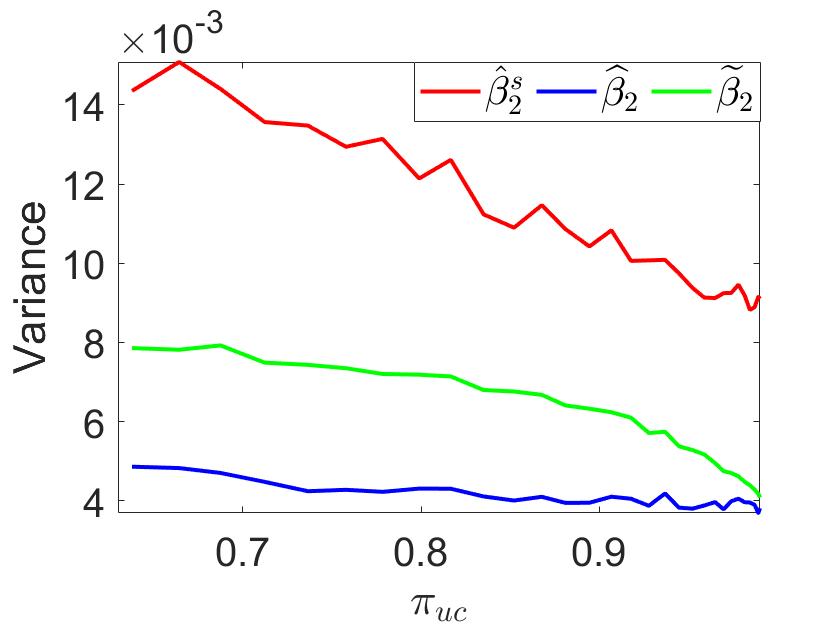}
  \caption{Variance of the estimators.}
      \label{fig:var}
      \end{minipage}
\end{figure}

\begin{figure}[H]
\begin{minipage}{0.46\textwidth}
 \centering
 \includegraphics[width=65mm]{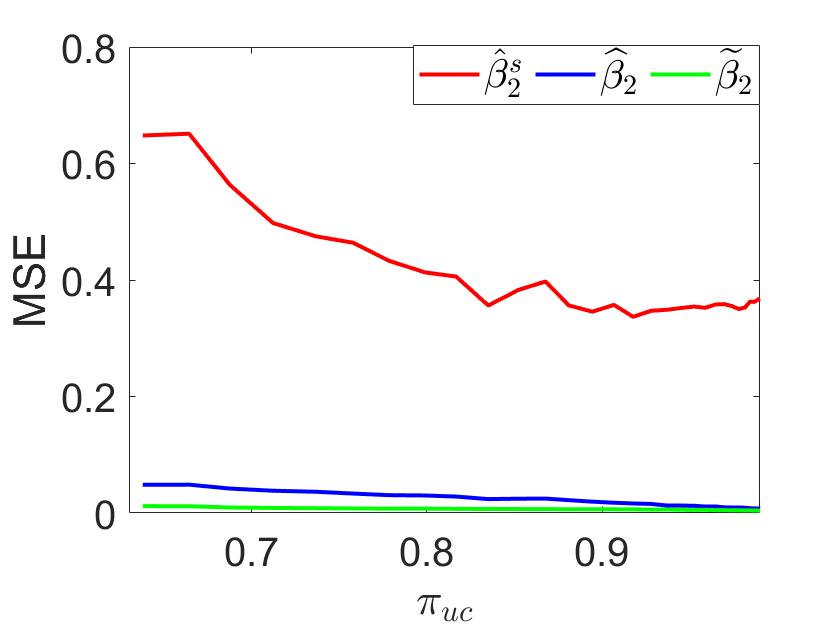}
\caption{MSE of the estimators.}
\label{fig:mse} 
\end{minipage}
  \quad\quad
  \begin{minipage}{0.46\textwidth}
  \centering
  \includegraphics[width=65mm]{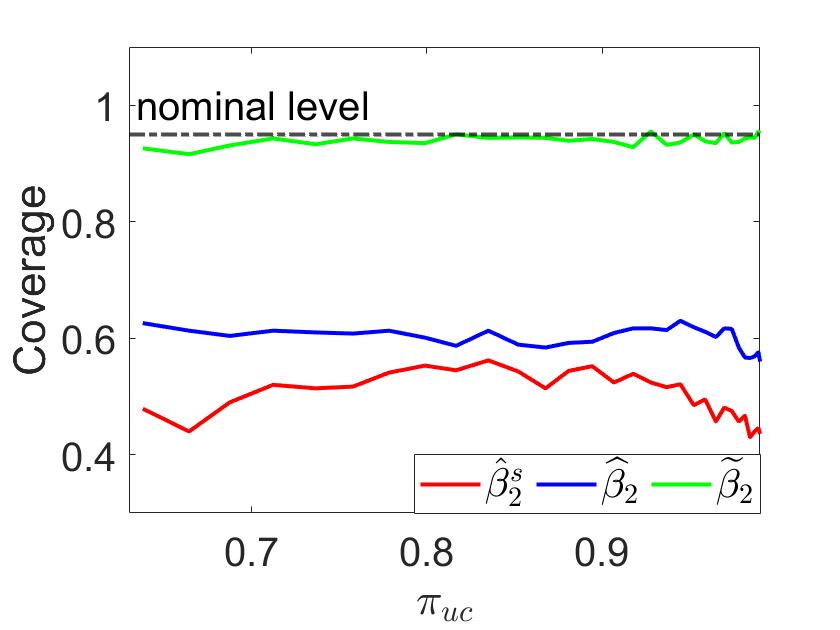}
   \caption{Coverage of 95\% confidence intervals.}
   \label{fig:cov}
   \end{minipage}
\end{figure}

We study $\widehat{\beta}$, $\widetilde{\beta}$ and the usual Stute estimator $$\widehat{\beta}^s= \argmin\limits_{b\in\R^p} \left|\left|Y^w-X^wb \right|\right|_2^2,$$ which is not robust to outliers. The algorithm of Section \ref{sec.algo} has been used to compute $\widehat{\beta}$, we initialize the search at ${a}^{(0)}=(0)_{1}^n$ and choose to stop after $10$ iterations, but the results are not sensitive to the number of iterations as long as it is large enough. The threshold of the two-step estimator $\widetilde{\beta}$ has been set at $\tau_0=0.3$. This decision was taken because $\widehat{\alpha}^w$ often has many small coefficients which are likely not outliers. For a grid of values of $\mu$ between $2$ and $5$ (with step $0.1$), we generate $1,000$ samples of size $1,000$ and compute these three estimators. The (estimated over all replications) probability $\pi_{uc}$ increases with $\mu$ from $64$ to $99\%$. 
We report the bias, variance and mean squared error (henceforth, MSE) of $\widehat{\beta}^s_2,\widehat{\beta}_2,\widetilde{\beta}_2$ in Figures \ref{fig:bias}, \ref{fig:var} and \ref{fig:mse}, respectively. In  Figure \ref{fig:cov}, we present the coverage of 95\% confidence intervals for $\beta_2$ based on the asymptotic normality of $\widehat\beta_2^s$, $\widehat\beta_2$ and $\widetilde{\beta}_2$ given in \cite{stute1996distributional}, Theorems \ref{an} and \ref{an2}, respectively. The estimator of the asymptotic variance is outlined in Section \ref{app.esas}. 

Let us analyze the results. As expected, Stute's estimator seems to be asymptotically biased. Our two-step estimator $\widetilde{\beta}_2$ has a better performance than $\widehat{\beta}_2$ in terms of bias, MSE and coverage. However, its variance is larger than that of $\widehat{\beta}_2$. Moreover, as anticipated, the precision of the three estimators deteriorates when the proportion of censored observations increases. Finally, the coverage of 95\% confidence intervals based on $\widetilde{\beta}_2$ is almost nominal, even for this relatively small sample size.

\appendix

\section{Technical conditions  and asymptotic variance}
\subsection{Conditions from \cite{stute1996distributional}}
\label{app.cond}

We focus on the case where the distribution of $\tilde Y_1$, that is $H(t)=\P(\tilde Y_1\le t)$, has no point mass at $\tau_{\tilde T_1}$, and treat $\xi$ as a covariate. In this setting the notations are less involved. For the general case see \cite{stute1996distributional}. We introduce the following quantities. For $t\in\R$, $x\in\R^p$ and $e\in\R$, let $F(x, e,  t) = \P( X_{11}\le x_1,\dots, X_{1p}\le x_p,\xi_1\le e , \tilde T_1 \le t)$, $G(t)=\P(C_1\le t)$ , $H^0(y)= \P(\tilde T_1\le t,\tilde \delta_1=0)$ and $H^{11}(x,e,t)=\P( X_{11}\le x_1,\dots, X_{1p}\le x_p, \xi_1\le e, \tilde T_1\le t, \tilde \delta_1=1)$. Moreover, for $y\in\R_+$, let $G(y-)$ be the left limit of $G$ at $y$.  Now, for $k=1,\dots, p$, and let also $\varphi_k:\ \R^p \times \R \mapsto \R : (x,e) \mapsto \varphi_k(x,e)=x_ke$. The conditions for asymptotic normality of the estimator of \cite{stute1996distributional} (in the absence of outliers) are as follows.

\begin{Assumption}\label{stute_rate} $H$ has no point mass at \red{$\tau_{\tilde T_1}$} and 
\begin{align}
\label{cond_stute_1}
& \int \int \int_0^{\tau_{\tilde T_1}}\varphi_k^2(x,e) \frac{dH^{11}(x,e,z)}{(1-G(z-))^2}<\infty;\\
& \notag \int \int \int_0^{\tau_{\tilde T_1}}  |\varphi_k(x,e)|\sqrt{\int_0^{z-} \frac{dG(y)}{(1-H(y))(1-G(y-))}} dF(x,e, z)<\infty,
\end{align}
for all $k=1,\dots, p$.
\end{Assumption}
Note that condition \eqref{cond_stute_1} is always satisfied when $\tau_{\red{\tilde T_1}}<\tau_{C_1}$.

\subsection{Asymptotic variance}\label{app.as}

For a function $\varphi:\R^p\times \R\mapsto \R$ and $t\in \R$, we define
\begin{align*}
\gamma_1^\varphi(t)&=\frac{1}{1-H(t)}\int \int \int 1_{\{t<z\}} \varphi(x,e) \frac{dH^{11}(x,e,z)}{1-G(z-)};\\
\gamma_2^\varphi(t)&=\int \int \int  \int 1_{\{z<t, z<y\}} \varphi(x,e) \frac{dH^0(z)dH^{11}(x,e,y)}{(1-H(z))^2(1-G(y-))}.
\end{align*}
For $k=1,\dots, p$, let $\psi_k= \varphi_k(X_1, \xi_1)(1-G(Y_1-))^{-1}\delta_1+ \gamma_1^{\varphi_k}(Y_1)(1-\delta_1) -\gamma_2^{\varphi_k}(Y_1)$. In \cite{stute1996distributional} it is shown that the asymptotic variance matrix of $\sqrt{n}\sum_{i=1}^n \tilde w_{(i)} \linebreak X_{(i)} \xi_{(i)}$ is given by  $\Sigma=(\Sigma_{k\ell})_{k,\ell=1}^p$, where $\Sigma_{k\ell}=\text{cov}(\psi_k,\psi_{\ell})$, and where $\tilde w_{(i)}$ is defined at the start of Appendix B.

\subsection{Estimation of the asymptotic variance}\label{app.esas}

The component $\Sigma_X$ can be naturally estimated by $(X^w)^\top X^w$. To estimate $\Sigma$, for $i=1,\dots, n$, we set $\widehat{\xi}_{(i)} = Y_{(i)} -X_{(i)}^\top\widehat{\beta} -\widehat{\alpha}_{(i)}$ for the first-step estimator and $\widetilde{\xi}_{(i)} = Y_{(i)} -X_{(i)}^\top\widetilde{\beta} -\widetilde{\alpha}_{(i)}$ for the second step estimator. Then, for $k=1,\dots,p $, we let
\begin{align*}
\widehat{\gamma}_1^{\varphi_k}(t)&=\frac{1}{1-\widehat{H}(t)}\frac{1}{n}\sum_{i=1}^n 1_{\{t<Y_{(i)} \}}  \frac{ \delta_{(i)}X_{(i)k}\widehat{\xi}_{(i)}}{1-\widehat{G}(Y_{(i)}-)};\\
\widehat{\gamma}_2^{\varphi_k}(t)&= \red{\frac{1}{n^2}}\sum_{i=1}^n \sum_{j=1}^n \int 1_{\{Y_{(j)}<t, Y_{(j)}<Y_{(i)}\}}   \frac{(1-\delta_{(j)})\delta_{(i)}X_{(i)k}\widehat{\xi}_{(i)}}{(1-\widehat{H}(Y_{(j)}))^2(1-\widehat{G}(Y_{(i)}-))},
\end{align*}
where $\widehat{G}$ is the Kaplan-Meier estimator of $G$ (using the observed sample) and $\widehat{H}(t)=n^{-1}\sum_{i=1}^n 1_{\{Y_i\le t\}}$. Also, we write $\widehat{\psi}_{ki}= X_{(i)k}\widehat{\xi}_{(i)}\delta_{(i)}(1-\widehat{G}(Y_{(i)}-))^{-1}+\widehat{ \gamma}_1^{\varphi_k}(Y_{(i)})(1-\delta_{(i)}) -\widehat{\gamma}_2^{\varphi_k}(Y_{(i)})$. Finally, the estimator $\widehat{\Sigma}$ of ${\Sigma}$ is the empirical covariance matrix of the $\widehat{\psi}_k,\ k=1,\dots, p$.

\section{Some results on the Kaplan-Meier weights}

Recall that $\tilde T_i=T_i-\alpha_i$, $\tilde\delta_i=1_{\{\tilde T_i\le C_i\}},\ \pi_{uc}=\P(\tilde \delta_1=1)$ and $\tilde Y_i =\min(\tilde T_i,C_i)$. Let $\tilde w_{(i)}$ be the Kaplan-Meier weight calculated based on the sample $\tilde S =\{(\tilde Y_i,\tilde \delta_i)\}_{i=1}^n$ of the observation whose rank equals $i$ in the sample $S=\{(Y_i,\delta_i)\}_{i=1}^n$. We have the following results.
\begin{Lemma}\label{new1} If $\pi_\alpha\to 0$, then $\widehat{\pi}_{uc} \xrightarrow{\P}\pi_{uc}$.
\end{Lemma}

\begin{Proof}
By the law of large numbers $n^{-1}\sum_{i=1}^n \tilde \delta_i\xrightarrow{\P}\pi_{uc}$. Moreover, $\sum_{i=1}^n\delta_i$ and $\sum_{i=1}^n \tilde \delta_i$ can only differ by at most $\norm{\alpha}_0$. Hence, we have $\widehat{\pi}_{uc}-n^{-1}\sum_{i=1}^n \tilde \delta_i=O_P(n^{-1}\norm{\alpha}_0)=o_P(1)$. 
\end{Proof}

\begin{Lemma} \label{new2}Assume that $\pi_\alpha\to 0$. For any $s<\pi_{uc}$ we have $ \norm{\tilde w}_{\infty}=o_P(n^{-s})$ and $ \norm{w}_{\infty}=o_P(n^{-s})$.
\end{Lemma}

\begin{Proof} The fact that $ \norm{\tilde w}_{\infty}=o_P(n^{-s})$ follows from Theorem 2.1 in \cite{shieh2000jump}. Then, by the same result, we have that for all $s<\pi_{ucn}$, $ \norm{w}_{\infty}=o_P(n^{-s})$, where $\pi_{ucn}=\P(\delta_i=1)$. To conclude the proof, note that $|\pi_{ucn}-\pi_{uc}|\le \pi_\alpha=o(1)$.
\end{Proof}

\begin{Lemma}\label{new3} Assume that $\pi_{\alpha}\to 0$. For any $s<\pi_{uc}$, we have $||w-\tilde w||_1=O_P(n^{1-s}\pi_\alpha (1+\log(n))).$
\end{Lemma}

\begin{Proof}
Let $\tilde Y_{\widetilde{(1)}},\dots, \tilde Y_{\widetilde{(n)}}$ be the order statistics of $\{\tilde Y_i\}_{i=1}^n$ and let $\tilde w_{\widetilde{(1)}},\dots, \tilde w_{\widetilde{(n)}}$ be the associated weights.  We now rewrite the expression of the weights as follows:
\begin{align*}  w_{(1)}= \frac{\delta_{(1)}}{n},\ w_{(i)}= \frac{\delta_{(i)}}{n}\prod_{j=1}^{i-1}\left( \frac{n-j+1}{n-j}\right)^{1-\delta_{(j)}},\  i=2,\dots,n;\\
\tilde w_{\widetilde{(1)}}= \frac{\tilde \delta_{\widetilde{(1)}}}{n},\ \tilde w_{\widetilde{(i)}}= \frac{\tilde \delta_{\widetilde{(i)}}}{n}\prod_{j=1}^{i-1}\left( \frac{n-j+1}{n-j}\right)^{1-\tilde \delta_{\widetilde{(j)}}},\  i=2,\dots,n.
\end{align*}
Assume now that there are $(n/2)-1\ge k\ge 1$ outliers. Let us consider an observation that is not an outlier and whose rank in sample $\{\tilde Y_i\}_{i=1}^n$ is $\tilde r\in\{k+1,\dots, n-k\}$.  Its rank $r$ in sample $\{Y_i\}_{i=1}^n$ then belongs to $\{\tilde r -k,\dots,\tilde r+k\} $.  We have
\begin{align*}w_{(r)}-\tilde w_{(r)}&=w_{(r)}-\tilde w_{\widetilde{(\tilde r)}}\\
&=  \frac{\delta_{(r)}}{n}\prod_{j=1}^{r-1}\left( \frac{n-j+1}{n-j}\right)^{1-\delta_{(j)}}-\frac{\tilde \delta_{\widetilde{(\tilde r)}}}{n}\prod_{j=1}^{\tilde r-1}\left( \frac{n-j+1}{n-j}\right)^{1-\tilde \delta_{\widetilde{(j)}}}\\
&=\frac{\tilde \delta_{\widetilde{(\tilde r)}}}{n} \left( \prod_{j=1}^{r-1}\left( \frac{n-j+1}{n-j}\right)^{1-\delta_{(j)}}-\prod_{j=1}^{\tilde r-1}\left( \frac{n-j+1}{n-j}\right)^{1-\tilde \delta_{\widetilde{(j)}}}\right).
\end{align*}
Since there are $k$ outliers, there are at most $k$ terms that differ between the products
$$\prod_{j=1}^{r-1}\left( \frac{n-j+1}{n-j}\right)^{1-\delta_{(j)}}\text{ and } \prod_{j=1}^{\tilde r-1}\left( \frac{n-j+1}{n-j}\right)^{1-\tilde \delta_{\widetilde{(j)}}}.$$
Let us denote by $A$ the common factors of both products, by $t_1,\dots, t_{k_1}$ the factors in the first product that are not in the second product, and by $\tilde t_1,\dots,\tilde t_{k_2}$ the factors in the second product that are not in the first product. We have that $k_1+k_2\le k$, and that
$$\prod_{j=1}^{r-1}\left( \frac{n-j+1}{n-j}\right)^{1-\delta_{(j)}}=A\prod_{i=1}^{k_1}t_{i} \text{ and } \prod_{j=1}^{\tilde r-1}\left( \frac{n-j+1}{n-j}\right)^{1-\tilde \delta_{\widetilde{(j)}}}= A\prod_{i=1}^{k_2}\tilde t_{i}.$$
Hence, it holds that 
\begin{align*} &\left|\prod_{j=1}^{r-1}\left( \frac{n-j+1}{n-j}\right)^{1-\delta_{(j)}}-\prod_{j=1}^{\tilde r-1}\left( \frac{n-j+1}{n-j}\right)^{1-\tilde \delta_{\widetilde{(j)}}}\right|\\
&=A\left|\prod_{i=1}^{k_1}t_{i} - \prod_{i=1}^{k_2}\tilde t_{i}\right|\\
& \le \prod_{j=1}^{\tilde r-1}\left( \frac{n-j+1}{n-j}\right)^{1-\tilde \delta_{\widetilde{(j)}}} \left(\max\left(\prod_{i=1}^{k_1}t_{i},\prod_{i=1}^{k_2}\tilde t_{i}\right)-1\right),
\end{align*}
where the last inequality is because $t_i\ge 1,\ i=1\dots,k_1$ and $\tilde t_i\ge 1,\ i=1,\dots,k_2$. 
Next, since 
$$j\mapsto \frac{n-j+1}{n-j}$$
is increasing in $j$, we have 
$$\max\left(\prod_{i=1}^{k_1}t_{i},\prod_{i=1}^{k_2}\tilde t_{i}\right) \le \prod_{j=\tilde r}^{\tilde r+k-1}\frac{n-j+1}{n-j} = \frac{n-\tilde r+1}{n-\tilde r-k+1}.$$
This yields 
\begin{align*}
& \prod_{j=1}^{r-1}\left( \frac{n-j+1}{n-j}\right)^{1-\delta_{(j)}}-\prod_{j=1}^{\tilde r-1}\left( \frac{n-j+1}{n-j}\right)^{1-\tilde \delta_{\widetilde{(j)}}} \\
& \le  \left(\frac{n-\tilde r+1}{n-\tilde r-k+1}-1\right)\prod_{j=1}^{\tilde r-1}\left( \frac{n-j+1}{n-j}\right)^{1-\tilde \delta_{\widetilde{(j)}}}. 
\end{align*}
We obtain that
$$|w_{(r)}-\tilde w_{(r)}| \le \tilde w_{(r)}\left(\frac{n-\tilde r+1}{n-\tilde r-k+1}-1\right)= \tilde w_{(r)}  \frac{k}{n-\tilde r-k+1}.$$
For the remaining observations (at most $2k$ with rank lower than $k$ or larger than $n-k$, plus $k$ outliers), we can bound $w_{(i)}-\tilde w_{(i)}$ by $\norm{w}_\infty+\norm{\tilde w}_{\infty}$. As a result, we have 
\begin{align*}\sum_{i=1}^n \left|\tilde w_{(i)}-w_{(i)}\right|&\le 3k \left(\norm{w}_\infty+\norm{\tilde w}_{\infty}\right)+ \norm{\tilde w}_{\infty}k\sum_{i=k+1}^{n-k} \frac{1}{n-i-k+1}\\
&\le 3k \left(\norm{w}_\infty+\norm{\tilde w}_{\infty}\right)+  \norm{\tilde w}_{\infty} k(1+\log(n)),
\end{align*}
where we use $$\sum_{i=k+1}^{n-k} \frac{1}{n-i-k+1}\le \sum_{i=1}^{n-2k}  \frac{1}{i}\le 1+\int_{1}^{n-2k-1} \frac1t dt\le (1+\log(n))$$
and the fact that $1/i\le \int_{i-1}^i dt/t $. 
To conclude the proof, use Lemma \ref{new2} and the fact that $\pi_\alpha\to 0$.\end{Proof}

\begin{Lemma}\label{new4} 
Let condition \eqref{Ei} hold. Then, $\sum_{i=1}^nw_{(i)}X_{(i)}X_{(i)}^\top\xrightarrow{\P}\Sigma_X .$
\end{Lemma}

\begin{Proof}By the Theorem in \cite{stute1993consistent}, we have that $\sum_{i=1}^n\tilde w_{(i)}X_{(i)}X_{(i)}^\top\xrightarrow{\P}\Sigma_X .$ Hence, it suffices to show that $\sum_{i=1}^n (w_{(i)}-\tilde w_{(i)})X_{(i)}X_{(i)}^\top=o_P(1)$.
This holds since $|\sum_{i=1}^n (w_{(i)}-\tilde w_{(i)})X_{(i)}X_{(i)}^\top|\le ||w-\tilde w||_1\norm{X}_{\infty}^2$, combined with Lemma \ref{new3} and condition \eqref{Ei}.
\end{Proof}

\begin{Lemma}\label{new5} Let Assumption \ref{DGP} (iv) hold. Then, $\sum_{i=1}^nw_{(i)}X_{(i)}\xi_{(i)}= \linebreak O_P(\max(n^{-1/2}, n^{1-\pi_{uc}+3\lambda_0}\norm{X}_\infty\pi_\alpha)).$
\end{Lemma}

\begin{Proof}By Theorem 1.1 in \cite{stute1996distributional}, we have $\sqrt{n}\sum_{i=1}^n\tilde w_{(i)}X_{(i)}\xi_{(i)}\xrightarrow{d}\mathcal{N}(0,\Sigma).$ Hence, it suffices to show that 
$\sum_{i=1}^n (w_{(i)}-\tilde w_{(i)})X_{(i)}\xi_{(i)}=O_P(n^{1-\pi_{uc}+3\lambda_0}\norm{X}_\infty \linebreak \pi_\alpha)$.
This follows from $|\sum_{i=1}^n (w_{(i)}-\tilde w_{(i)})X_{(i)}\xi_{(i)}| \le ||w-\tilde w||_1\norm{X}_{\infty}\norm{\xi}_\infty$, together with Lemma \ref{new3} and Assumption \ref{DGP} (iv).
\end{Proof}

\begin{Lemma}\label{new6} Let Assumption \ref{DGP} (iv) and condition \eqref{Eii} hold.  Then, $$\sqrt{n}\sum_{i=1}^n w_{(i)} X_{(i)}\xi_{(i)}\xrightarrow{d} \mathcal{N}(0, \Sigma).$$
\end{Lemma}

\begin{Proof}Condition \eqref{Eii} implies that $\left|\sum_{i=1}^n (w_{(i)}-\tilde w_{(i)})X_{(i)}\xi_{(i)}\right|=o_P(n^{-1/2}).$
This together with Theorem 1.1 in \cite{stute1996distributional} gives the result.
\end{Proof}

\section{Proofs} 
\label{sec.app}

In this section, we prove the different theorems of the paper. The proofs rely on results on the penalty level $\lambda$ and the estimation error $\widehat{\alpha}-\alpha$ which are given in the next two subsections.

\subsection{Penalty level} 

Let us first prove a number of lemmas.

\begin{Lemma}
\label{penalty}
Under Assumption \ref{DGP}, \ref{effective} and when $\pi_\alpha\to 0$,  we have $\lambda= o_P(n^{2\lambda_0-(\pi_{uc}/2)})$.
\end{Lemma}
\begin{Proof}
By Lemma \ref{new1}, $\widehat{\pi}_{uc}\xrightarrow{\P} \pi_{uc}$, which implies that  
$$\lim_{n\to\infty}\P(\pi_{uc}-\widehat{\pi}_{uc}-2\lambda_0<-\lambda_0)=1.$$ This leads to $\lambda=n^{\lambda_0-\widehat{\pi}_{uc}/2}=n^{2\lambda_0-(\pi_{uc}/2)}n^{(\pi_{uc}-\widehat{\pi}_{uc})/2-\lambda_0}=o_P(n^{2\lambda_0-(\pi_{uc}/2)})$.
\end{Proof}

\begin{Lemma}
\label{choice}
Under Assumption \ref{DGP} and condition \eqref{Ei},  it holds that \linebreak $\norm{P_{X^w}\xi^w}_\infty= O_P(n^{-1/2})$.
\end{Lemma}

\begin{Proof}
By Lemma \ref{new4}, we have $(X^w)^\top X^w=\sum_{i=1}^nw_{(i)}X_{(i)}X_{(i)}^\top\xrightarrow{\P} \Sigma_X.$
This implies that $\norm{X^w}_2=O_P(1)$. Moreover, by the continuous mapping theorem, we have $ \left((X^w)^\top X^w\right)^{-1}\xrightarrow{\P} \Sigma_X^{-1}$. Moreover Lemma \ref{new5} and the fact that condition \eqref{Ei} holds yield $\sqrt{n}(X^w)^\top\xi^w =O_P(1)$. This leads to
$$\sqrt{n}\norm{P_{X^w}\xi^w}_2\le\norm{X^w}_2\norm{((X^w)^\top X^w)^{-1}}_2\norm{\sqrt{n}(X^w)^\top\xi^w}_2=O_P(1),$$
which concludes the proof because the sup-norm is smaller than the $\ell_2$-norm.
\end{Proof}

\begin{Lemma}
\label{correffective}
Under Assumptions \ref{DGP}, \ref{effective}  and condition \eqref{Ei},  it holds that $\lim\limits_{n\to \infty}\P(\lambda\ge 4\norm{M_{X^w}\xi^w}_\infty)=1$. 
\end{Lemma}

\begin{Proof}
By Lemma \ref{choice} and the fact that $\norm{\xi}_{\infty}=O_P(n^{\lambda_0})$ (Assumption \ref{DGP} (iv)), we have
\begin{align} \notag 4\norm{M_{X^w}\xi^w}_{\infty}&=O_P(\norm{\xi^w}_{\infty}+\norm{P_{X^w}\xi^w}_\infty)\\*
\notag &= O_P(\sqrt{\norm{w}_\infty}\norm{\xi}_\infty)+O_P(n^{-1/2}) \\
\label{fbound}&=O_P(\sqrt{\norm{w}_\infty}n^{\lambda_0})+O_P(n^{-1/2}).\end{align}
Next, take $s<\pi_{uc}-\lambda_0$. By Theorem 2.1 in  \cite{shieh2000jump}, it holds that
\begin{equation} \label{sbound}  \norm{w}_\infty = o_P(n^{-s}) = o_P(n^{2\lambda_0-\widehat{\pi}_{uc}}n^{-s+\widehat{\pi}_{uc}-\lambda_0}). \end{equation}
Since by the law of large numbers $\widehat{\pi}_{uc}\xrightarrow{\P} \pi_{uc}$, we have $\lim\limits_{n\to\infty}\P(s-\widehat{\pi}_{uc}+\lambda_0<0)=1$ and, therefore, $n^{-s+\widehat{\pi}_{uc}-\lambda_0}=o_P(1)$. This, \eqref{fbound} and \eqref{sbound} imply $4\norm{M_{X^w}\xi^w}_{\infty}=o_P(n^{ \lambda_0-(\widehat{\pi}_{uc}/2)})=o_P(\lambda)$.
\end{Proof}

\subsection{Bound on the estimation error of $\alpha$}

In this subsection, we derive a bound on $\norm{\widehat{\alpha}^w-\alpha^w}_1$. Define $J=\{i\in\{1,\dots,n\}:\ |\alpha_{(i)}|>0\}$ and the following cone
$$C=\left\{\delta\in\R^n:\ \norm{\delta_{J^c}^w}_1\le \frac{3}{2} \norm{\delta_{J}^w}_1 \right\}.$$
We introduce the following compatibility constant:
$$\kappa=\inf_{\delta\in C}\frac{\norm{M_{X^w}\delta^w }_2}{\norm{\delta^w}_2}.$$
Let us show the following lemmas.

\begin{Lemma}\label{compatibility}
Let Assumptions \ref{DGP} (iii), \ref{censoring} and condition (\ref{Ei}) hold. Then, \linebreak $\lim\limits_{n\to\infty}\P(\kappa>1/2)=1$.
\end{Lemma}

\begin{Proof} Note that by Lemma \ref{new4} we have that $(X^w)^\top X^w$ converges in probability to $\Sigma_X$.  Since $\Sigma_X$ is positive definite and the determinant of a matrix is a continuous mapping, $\lim_{n\to\infty}\P(\mathcal{X})=1$, where $\mathcal{X}$ is the event that $(X^w)^\top X^w$ is positive definite. In the remainder of this proof, let us work on the event $\mathcal{X}$. We have
$M_{X^w}\delta^w=\delta^w -X^w((X^w)^\top X^w)^{-1}(X^w)^{\top}\delta^w,$ which implies that
\begin{align}
\notag\norm{M_{X^w}\delta^w}_2&\ge \norm{\delta^w}_2-\norm{X^w((X^w)^\top X^w)^{-1}(X^w)^{\top}\delta^w}_2\\
\notag &= \norm{\delta^w}_2-\left|\left|\sum_{k=1}^{p}X_{k\cdot }^w\left(((X^w)^\top X^w)^{-1}(X^w)^{\top}\delta^w\right)_k\right\vert\right\vert_2\\
\notag &\ge \norm{\delta^w}_2-\sum_{k=1}^{p}\left|\left|X_{k\cdot }^w\left(((X^w)^\top X^w)^{-1}(X^w)^{\top}\delta^w\right)_k\right\vert\right\vert_2\\
\notag &\ge \norm{\delta^w}_2-\sum_{k=1}^{p}\left|\left|X_{k\cdot }^w\right|\right|_2\left|\left|((X^w)^\top X^w)^{-1}(X^w)^{\top}\delta^w\right\vert\right\vert_{\infty}\\
\notag &\ge \norm{\delta^w}_2-\sum_{k=1}^{p}\left|\left|X_{k\cdot }^w\right|\right|_2\left|\left|((X^w)^\top X^w)^{-1}(X^w)^{\top}\delta^w\right\vert\right\vert_{2}\\
\notag &\ge \norm{\delta^w}_2-\sum_{k=1}^{p}\norm{X_{k\cdot }^w}_{2}\left\vert \left\vert ((X^w)^\top X^w)^{-1}\right\vert \right\vert_{\text{op}}\norm{(X^w)^{\top}\delta^w}_2.
\end{align}
Note that 
\begin{align*}
\norm{(X^w)^{\top}\delta^w}_2 &= \sqrt{\sum_{k=1}^p((X_{k\cdot }^w)^{\top}\delta^w)^2} \\
& =\sqrt{\sum_{k=1}^p\left(\sum_{i=1}^n X^w_{ki}\delta^w_i\right)^2} 
\le \sqrt{p}\norm{X^w}_\infty \norm{\delta^w}_1.
\end{align*}
This leads to
\begin{align}
& \notag \norm{M_{X^w}\delta^w}_2 \\
&\ge\norm{\delta^w}_2-\sum_{k=1}^{p}\norm{X_{k\cdot }^w}_{2}\left\vert \left\vert ((X^w)^\top X^w)^{-1}\right\vert \right\vert_{\text{op}} \sqrt{p}\norm{X^w}_\infty \norm{\delta^w}_1 \\
\notag &\ge \norm{\delta^w}_2-\sum_{k=1}^{p}\norm{X_{k\cdot }^w}_{2}\left\vert \left\vert ((X^w)^\top X^w)^{-1}\right\vert \right\vert_{\text{op}} \sqrt{p}\norm{X^w}_\infty \frac52\norm{\delta_J^w}_1 \\
\notag   &\ge \norm{\delta^w}_2-\sum_{k=1}^{p}\norm{X_{k\cdot }^w}_{2}\left\vert \left\vert((X^w)^\top X^w)^{-1}\right\vert \right\vert_{\text{op}} \sqrt{p}\norm{X^w}_\infty \frac52\sqrt{\norm{\alpha}_0}\norm{\delta_J^w}_2,\end{align}
where the second inequality is due to the fact that $\delta\in C$.
Then, we have, by definition of $\kappa$ and the fact that $\norm{\delta^w}_2\ge \norm{\delta^w_J}_2$,
\begin{equation}\label{lastcompatibility}
 \kappa 
\ge 1-\sum_{k=1}^{p}\norm{X_{k\cdot }^w}_{2}\left\vert \left\vert \left((X^w)^{\top}(X^w)\right)^{-1}\right\vert \right\vert_{\text{op}} \frac52\sqrt{{p}}\sqrt{\frac{\norm{\alpha}_0}{n}}\sqrt{n}\sqrt{\norm{w}_\infty}\norm{X}_{\infty}.
\end{equation}
 Now, since $(X^w)^\top X^w$ converges in probability to $\Sigma_X$, we have $$\left\vert \left\vert \left( (X^w)^{\top}X^w\right)^{-1}\right\vert \right\vert_{\text{op}}=O_{P}(1)$$ and $$\sum_{k=1}^{p}\norm{X_{k\cdot }^w}_2=\sum_{k=1}^{p}\sqrt{\left((X^w)^{\top}X^w\right)_{kk}}=O_{P}(1),$$
both implying that $\sum_{k=1}^{p}\norm{X_{k\cdot }^w}_2\vert \vert ((X^w)^{\top}X^w)^{-1}\vert \vert_{\text{op}}=O_{P}(1).$
 By Theorem 2.1 in  \cite{shieh2000jump}, it holds that $\norm{w}_{\infty}=o_P(n^{-(\pi_{uc}-\lambda_0)})$. Then, we obtain 
\begin{align*}&\sum_{k=1}^{p}\norm{X_{k\cdot }^w}_{2}\left\vert \left\vert \left((X^w)^{\top}(X^w)\right)^{-1}\right\vert \right\vert_{\text{op}} \frac52\sqrt{{p}}\sqrt{\frac{\norm{\alpha}_0}{n}}\sqrt{n}\sqrt{\norm{w}_\infty}\norm{X}_{\infty}\\
&=O_P\left(n^{(1-\pi_{uc}+\lambda_0)/2}\sqrt{\pi_\alpha}\norm{X}_\infty\right)=o_P(1),\end{align*}
by condition (\ref{Ei}). This yields the result by \eqref{lastcompatibility}.
\end{Proof}

\begin{Lemma}\label{consistent}
Let Assumptions \ref{DGP}, \ref{effective}, \ref{censoring} and condition \eqref{Ei} hold.  Then, 
$$\norm{\widehat{\alpha}^w-\alpha^w}_1=O_P(n\pi_\alpha\lambda).$$
\end{Lemma}

\begin{Proof}
We work on the event 
$$\mathcal{E}=\Big\{\lambda \ge 4\norm{M_{X^w}\xi^w}_{\infty}\Big\}\cup \left\{ \kappa \ge\frac12\right\},$$
whose probability goes to $1$ by Lemmas \ref{correffective} and \ref{compatibility}.
Let us define $\Delta=\widehat{\alpha}-\alpha$. Now, remark that
\begin{align}
 \notag \norm{\alpha^w}_1- \norm{\widehat{\alpha}^w}_1&= \norm{\alpha^w}_1-\norm{\alpha^w+\Delta^w}_1\\
\notag&= \norm{\alpha^w}_1-\norm{\alpha^w+\Delta_J^w}_1-\norm{\Delta_{J^c}^w}_1\\
\label{boundalpha} &\le \norm{\Delta_J^w}_1-\norm{\Delta_{J^c}^w}_1.
\end{align}
By \eqref{progalpha}, we have 
\begin{equation} \label{lbound}\norm{M_{X^w}(Y^w-\widehat{\alpha}^w)}_{2}^2- \norm{M_{X^w}(Y^w-\alpha^w)}_{2}^2  \le  \lambda (\norm{\alpha^w}_{1}-\norm{\widehat{\alpha}^w}_1). \end{equation}
 By convexity of $a\in\R^n\mapsto \norm{M_{X^w}(Y^w-a^w)}_{2}^2$, it holds that
\begin{align}\notag  \norm{M_{X^w}(Y^w-\widehat{\alpha}^w)}_{2}^2- \norm{M_{X^w}(Y^w-\alpha^w)}_{2}^2&\ge -2\left<M_{X^w}\xi^w,\Delta^w\right>\\
&\ge -2\norm{M_{X^w}\xi^w}_\infty\norm{\Delta^w}_1 \notag\\
\label{rbound}&\ge  -\frac{\lambda}{2} \norm{\Delta^w}_1,
\end{align}
where the last inequality is thanks to Lemma \ref{correffective}.
Combining \eqref{boundalpha}, \eqref{lbound} and \eqref{rbound}, we get $\norm{\Delta_{J^c}^w}_1\le 3\norm{\Delta_{J}^w}_1$
which implies that $\Delta \in C$.
Next, by H\"older's inequality, we have 
   \begin{align*} \norm{M_{X^w}(Y^w-\widehat{\alpha}^w)}_{2}^2- \norm{M_{X^w}(Y^w-\alpha^w)}_{2}^2&=\norm{M_{X^w}\Delta^w}_{2}^2-2\left<M_{X^w}\xi^w,\Delta^w\right>\\
 &\ge \norm{M_{X^w}\Delta^w}_{2}^2-\frac{\lambda}{2} \norm{\Delta^w}_1.
\end{align*}
Combining this, \eqref{boundalpha} and \eqref{lbound}, we obtain $(3/2)\lambda\norm{\Delta_J^w}_1\ge \norm{M_{X^w}\Delta^w}_{2}^2.$
Therefore, we have $(3/2) \lambda\norm{\Delta_J^w}_1\ge (1/4)\norm{\Delta^w}_2^2$ by the definition of $\kappa$, the fact that $\Delta\in C$ and Lemma \ref{compatibility}. Since $\norm{\Delta_J^w}_1\le \sqrt{\norm{\alpha}_{0}}\norm{\Delta_J^w}_2\le \sqrt{\norm{\alpha}_{0}}\norm{\Delta^w}_2$, this implies $\norm{\Delta^w}_2\le 6\lambda \sqrt{\norm{\alpha}_0}=O_P(\sqrt{n\pi_{\alpha}}\lambda)$ and, then, $\norm{\Delta^w}_1\le 4\norm{\Delta_J}_1\le \sqrt{\norm{\alpha}_{0}}\norm{\Delta^w}_2  =O_P(n\pi_{\alpha}\lambda).$
\end{Proof}

\subsection{Proof of Theorems \ref{consth} and \ref{an}}\label{sec.proof_an}

First, we prove Theorem \ref{consth}. Using \eqref{wOLS}, we have 
\begin{align} \label{decompn}
\widehat{\beta}-\beta = & \left((X^w)^\top X^w\right)^{-1}\sum_{i=1}^n w_{(i)} X_{(i)}\xi_{(i)} \\
& + \left((X^w)^\top X^w\right)^{-1}\sum_{i=1}^n w_{(i)} X_{(i)}(\alpha_{(i)}-\widehat{\alpha}_{(i)}). \nonumber \end{align}
By Lemma \ref{new4} and the continuous mapping theorem, $ \left((X^w)^\top X^w\right)^{-1}\xrightarrow{\P} \Sigma_X^{-1}$. By Lemma \ref{new5}, it holds that 
\begin{equation}\label{ann} \sum_{i=1}^nw_{(i)}X_{(i)}\xi_{(i)}=O_P\left(\max\left(\frac{1}{\sqrt{n}}, n^{1-\pi_{uc}+3\lambda_0}\norm{X}_\infty\pi_\alpha\right)\right).\end{equation}
Next, we have 
\begin{align*}\left|  \sum_{i=1}^n w_{(i)} X_{(i)}(\alpha_{(i)}-\widehat{\alpha}_{(i)})\right|&\le\sum_{k=1}^p\left|\sum_{i=1}^n(X_{ki }^w)^{\top}\sqrt{w_{(i)}} (\alpha_{(i)}-\widehat{\alpha}_{(i)})\right|
\\
&\le p\sqrt{\norm{w}_{\infty}}\norm{X}_ \infty \norm{\widehat{\alpha}^w-\alpha^w}_1.
\end{align*}
From Lemma \ref{new2} it follows that $\sqrt{\norm{w}_{\infty}}=o_P(n^{-(\pi_{uc}-2\lambda_0)/2})$. From Lemmas \ref{penalty} and \ref{consistent}, we obtain 
\begin{align}
\notag  \left|  \sum_{i=1}^n w_{(i)} X_{(i)}(\alpha_{(i)}-\widehat{\alpha}_{(i)})\right|&=O_P(n^{(-\pi_{uc}+2\lambda_0)/2}\norm{X}_{\infty}n\lambda\pi_\alpha)\\
\label{neg}&=O_P(n^{1-\pi_{uc}+3\lambda_0}\norm{X}_\infty\pi_\alpha).
\end{align} 
This, \eqref{decompn}, \eqref{ann} and the fact that $ \left((X^w)^\top X^w\right)^{-1}\xrightarrow{\P} \Sigma_X^{-1}$ show Theorem \ref{consth}.

Let us now prove Theorem \ref{an}. By Lemma \ref{new6} we have $\sqrt{n}\sum_{i=1}^n w_{(i)} X_{(i)}\xi_{(i)} \linebreak \xrightarrow{d} \mathcal{N}(0, \Sigma)$.
By condition (\ref{Eii}),  \eqref{neg} implies that $$\sqrt{n} \left|  \sum_{i=1}^n w_{(i)} X_{(i)}(\alpha_{(i)}-\widehat{\alpha}_{(i)})\right|=o_P(1).$$
This, \eqref{decompn} and Slutsky's theorem prove that $\sqrt{n} (\widehat{\beta}-\beta )\xrightarrow{d}  \mathcal{N}(0, \Sigma_X^{-1}\Sigma\Sigma_X^{-1})$.

\subsection{Proof of Theorem \ref{an2}}

Let us introduce $\widehat{J}=\widehat{J}(0)$. Remark that for all $\tau_0\ge 0$, $\widehat{J}(\tau_0)\subset\widehat{J}$, which implies that all the nonzero coefficients of $\widetilde{\alpha}_{\widehat{J}(\tau_0)}^w-\widehat{\alpha}^w$ belong to $\widehat{J}$. The proof proceeds in several parts.

\subsubsection{Sparsity bound}

First, we establish a bound on the number of nonzero coefficients of $\widehat{\alpha}$.

\begin{Lemma}\label{sparsity_bound}
Let Assumptions \ref{DGP}, \ref{effective}, \ref{censoring} and condition \eqref{Ei} hold. Then, 
$$|\widehat{J}|=O_P(n\pi_{\alpha}).$$
\end{Lemma}

\begin{Proof}
 By the Karush-Kuhn-Tucker conditions of program \eqref{progalpha}, we have $$2\left|\left(M_{X^w}(Y^w-\widehat{\alpha}^w)\right)_{(i)}\right|=\lambda$$
for all $i\in\widehat{J}$.
Summing over $i\in \widehat{J}$, we get
\begin{align*}
\sqrt{|\widehat{J}|}\lambda &=2\sqrt{\sum_{i\in \widehat{J}}\left(M_{X^w}(Y^w-\widehat{\alpha}^w)\right)_{(i)}^2} \\
&\le 2\sqrt{\sum_{i\in\widehat{J}
}\left(M_{X^w}\xi^w\right)_{(i)}^2}+2 \sqrt{\sum_{i\in \widehat{J}}\left(M_{X^w}(\widehat{\alpha}^w-\alpha^w)\right)_{(i)}^2}\\
&\le 2 \norm{M_{X^w}\xi^w}_\infty \sqrt{|\widehat{J}|}+2 \norm{M_{X^w}(\widehat{\alpha}^w-\alpha^w)}_2\\
&\le \sqrt{|\widehat{J}|}O_P\left(\frac{\lambda}{2}\right) + O_P\left(\lambda\sqrt{n\pi_\alpha}\right),
\end{align*}
where the first inequality is due to the triangular inequality and the last inequality follows from Assumption \ref{effective} and Lemmas  \ref{correffective} and \ref{consistent}.
\end{Proof}

\subsubsection{Penalty level} 

We prove that the penalty level is large compared to $\norm{M_{X^w}\widehat{\xi}^w
w}_\infty$.

\begin{Lemma}\label{correffective2}
Under Assumptions \ref{DGP} and \ref{effective}, we have 
$$\lim\limits_{n\to \infty}\P(\lambda\ge (3/2)\norm{M_{X^w}\widehat{\xi}^w}_\infty)=1. $$
\end{Lemma}

\begin{Proof}
By the Karush-Kuhn-Tucker conditions of program \eqref{progalpha}, we have that $2M_{X^w}(\alpha^w-\widehat{\alpha}^w)=\lambda\widehat{z}$, where $\widehat{z}\in \partial |\widehat{\alpha}^w|_1$, the subdifferntial of the $\ell_1$-norm at the point $\widehat{\alpha}^w$. Since $|\widehat{z}|_\infty \le 1$ (see Lemma D.5 in \cite{giraud2014introduction}), we obtain $\norm{M_{X^w}(\alpha^w-\widehat{\alpha}^w)}_\infty \le \lambda/2$. Hence, it holds that 
\begin{align*}
\norm{M_{X^w}\widehat{\xi}^w}_\infty&\le \norm{M_{X^w}\xi^w}_\infty+\norm{M_{X^w}(\alpha^w-\widehat{\alpha}^w)}_\infty\\
&\le \norm{M_{X^w}\xi^w}_\infty+\frac{\lambda}{2}.
\end{align*}
This implies that $\lambda-(3/2)\norm{M_{X^w}\widehat{\xi}^w}_\infty\ge (1/4)\lambda-\norm{M_{X^w}\xi^w}_\infty$, where the probability that the right-hand side is positive goes to $1$ by Lemma \ref{correffective}.
\end{Proof}

\subsubsection{Sparse eigenvalues}

For $s\in \N_*$, we define the following quantity which measures the effect of applying $M_{X^w}$ to a sparse vector:
$$\kappa_s=\inf_{\delta\in \R^p:\ \norm{\delta}_0\le s}\frac{\norm{M_{X^w}\delta^w }_2}{\norm{\delta^w}_2}.$$
We have the following result on $\kappa_{|\widehat{J}|}$.

\begin{Lemma}\label{sparse_eng}
Let Assumptions \ref{DGP}(iii), \ref{censoring} and condition (\ref{Ei}) hold.  Then, $$\lim_{n\to\infty}\P\Big(\kappa_{|\widehat{J}|}>1/2\Big)=1.$$
\end{Lemma}

\begin{Proof} Take $\delta\in \R^p$ such that $\norm{\delta}_0\le |\widehat{J}|$. Following the beginning of the proof of Lemma \ref{compatibility}, we have
\begin{align}
\notag \norm{M_{X^w}\delta^w}_2 &\ge\norm{\delta^w}_2-\sum_{k=1}^{p}\norm{X_{k\cdot }^w}_{2}\left\vert \left\vert ((X^w)^\top X^w)^{-1}\right\vert \right\vert_{\text{op}} \sqrt{p}\norm{X^w}_\infty \norm{\delta^w}_1 \\
\notag &\ge \norm{\delta^w}_2-\sum_{k=1}^{p}\norm{X_{k\cdot }^w}_{2}\left\vert \left\vert ((X^w)^\top X^w)^{-1}\right\vert \right\vert_{\text{op}} \sqrt{p}\norm{X^w}_\infty \sqrt{|\widehat{J}|}\norm{\delta^w}_2.\end{align}
By definition of $\kappa_{|\widehat{J}|}$, we obtain
$$
\kappa_{|\widehat{J}|}
\ge 1-\sum_{k=1}^{p}\norm{X_{k\cdot }^w}_{2}\left\vert \left\vert \left((X^w)^{\top}X^w\right)^{-1}\right\vert \right\vert_{\text{op}} \sqrt{{p}}\sqrt{\frac{|\widehat{J}|}{n}}\sqrt{n}\sqrt{\norm{w}_\infty}\norm{X}_{\infty}.
$$
We conclude as in the proof of Lemma \ref{compatibility} (using Lemma \ref{sparsity_bound} to bound $|\widehat{J}|$).
\end{Proof}

\subsubsection{Rate of convergence of $\boldsymbol{\widetilde{\alpha}_{\widehat{J}(\tau_0)}^w}$}

We have the following lemma which bounds the difference between $\widetilde{\alpha}_{\widehat{J}(\tau_0)}^w$ and $\widehat{\alpha}^w$. 

\begin{Lemma}\label{consistent2}
Let Assumptions \ref{DGP}, \ref{effective}, \ref{censoring} and condition (\ref{Ei}) hold.  Then, 
$$\norm{\widetilde{\alpha}_{\widehat{J}(\tau_0)}^w-\widehat{\alpha}^w}_1=O_P(n\pi_\alpha\lambda).$$
\end{Lemma}

\begin{Proof}
We work on the event 
$$\mathcal{E}=\left\{\lambda \ge \frac32\norm{M_{X^w}\xi^w}_{\infty}\right\}\cup \left\{ \kappa_{|\widehat{J}|} \ge\frac12\right\},$$
whose probability goes to $1$, by Lemmas \ref{correffective2} and \ref{sparse_eng}.
Let $\Delta^w=\widetilde{\alpha}_{\widehat{J}}^w-\widehat{\alpha}^w$. Using \eqref{equiv_formulation}, we have
   \begin{align*} 0&\ge \norm{M_{X^w}(Y^w-\widetilde{\alpha}_{\widehat{J}(\tau_0)}^w)}_{2}^2- \norm{M_{X^w}(Y^w-\widehat{\alpha}^w)}_{2}^2\\
   &=\norm{M_{X^w}\Delta^w}_{2}^2-2\left<M_{X^w}\widehat{\xi}^w,\Delta^w\right>\\
 &\ge \frac{1}{4} \norm{\Delta^w}_2^2-3\lambda \norm{\Delta^w}_1\\
&\ge \frac{1}{4} \norm{\Delta^w}_2^2-3\sqrt{|\widehat{J}|} \norm{\Delta^w}_2.
\end{align*}
This yields $ \norm{\Delta^w}_2\le 12\lambda \sqrt{|\widehat{J}|}=O_P(\sqrt{n\pi_\alpha}\lambda)$ by Lemma \ref{sparsity_bound}, which implies $\norm{\Delta^w}_1\le \sqrt{|\widehat{J}|}\norm{\Delta^w}_2=O_P(n\pi_\alpha\lambda)$ by Lemma \ref{sparsity_bound}. 
\end{Proof}

\subsubsection{End of the proof of Theorem \ref{an2}}

By \eqref{wOLS} and \eqref{equiv_formulation}, we have
$$\widetilde{\beta}-\widehat{\beta} = \left((X^w)^\top X^w\right)^{-1}\sum_{i=1}^n w_{(i)} X_{(i)}(\widehat{\alpha}_{(i)}-\widetilde{\alpha}_{(i)}).$$
We obtain the result following the arguments of Section \ref{sec.proof_an} and using Lemma \ref{consistent2}.

% Main text entry area

 \begin{acks}[Acknowledgments]
 Financial support from the European Research Council (2016-2021, Horizon 2020 / ERC grant agreement No.\ 694409) is gratefully acknowledged.
\end{acks}

%%%%%%%%%%%%%%%%%%%%%%%%%%%%%%%%%%%%%%%%%%%%%%
%% Supplementary Material, if any, should   %%
%% be provided in {supplement} environment  %%
%% with title and short description.        %%
%%%%%%%%%%%%%%%%%%%%%%%%%%%%%%%%%%%%%%%%%%%%%%
%\begin{supplement}
%\stitle{???}
%\sdescription{???.}
%\end{supplement}

%% if your bibliography is in bibtex format, uncomment commands:
\bibliographystyle{imsart-number} % Style BST file (imsart-number.bst or imsart-nameyear.bst)
\bibliography{ref}       % Bibliography file (usually '*.bib')

%% or include bibliography directly:
% \begin{thebibliography}{}
% \bibitem{b1}
% \end{thebibliography}

\end{document}